\documentclass[reqno]{article}
\usepackage{ae} 
\usepackage[T1]{fontenc}
\usepackage[ansinew]{inputenc}
\usepackage{amsmath}
\usepackage{graphicx}
\usepackage{color}
\usepackage[colorlinks]{hyperref}

\newtheorem{thm}{Theorem}
\newtheorem{lem}[thm]{Lemma}

\newtheorem{rmk}[thm]{Remark}

\begin{document}

\centerline{\bf Van der Waerden's Theorem on Homothetic Copies of
$\{1,1+s, 1+s+t \}$}

\medskip
\centerline{{\bf Byeong Moon Kim \footnote{Department of
Mathematics,
 Kangnung National University,
 Kangnung, Kangweondo, KOREA \\E-mail: bmk@knusun.kangnung.ac.kr } and Yoomi Rho }\footnote{
Department of Mathematics, Incheon University, Incheon
402-749, Korea\\
 E-mail: rho@incheon.ac.kr }}

\begin{abstract}
 For all positive
integers $s$ and $t$, Brown et. al \cite{B} defined $f(s,t)$ to be
the smallest positive integer $N$ such that every $2$-coloring of
$[1, N]$ has a monochromatic homothetic copy of $\{1,1+s, 1+s+t
\}$. They proved that $f(s,t) \leq 4(s+t)+1$ for all $s$, $t$ and
that the equality holds in the case where both
$s/g\not\equiv0\pmod 4$ and $t/g\not\equiv0\pmod 4$ with
$g=gcd(s,t)$
 and in
many other cases. Also they proved that for all positive integer
$m$, $f(4mt,t)=f(t,4mt)=4(4mt+t)-t+1$ or $4(4mt+t)+1$. In this
paper, we show that $f(4mt,t)=f(t,4mt)=4(4mt+t)-t+1$ and that for
all the other $(s,t)$, $f(s,t)=4(s+t)+1$.
\end{abstract}

keywords: van der Waerden's theorem, arithmetic progression,
homothetic copy, $2-$coloring, monochromatic triple

\baselineskip 0.3in
\section{Introduction}

 Van der Waerden's
theorem on arithmetic progressions \cite{W} states that for every
positive integer $k$ there is a smallest positive integer $w(k)$
such that every $k$-coloring of $[1,w(k)]=\{1,2,\ldots, w(k)\}$
has a monochromatic $k$-term arithmetic progression. There are
lots of results on the estimation of $w(k)$ for large $k$.( See
\cite{GRS} ]. )

 Let $s$, $t$, $m$ be positive integers. A homothetic
copy of $\{1,1+s, 1+s+t \}$ is any set of the form $\{x,x+ys,
x+ys+yt \}$ where $x$ and $y$ are positive integers. 
Regarding $3$-term arithmetic progressions as homothetic copies of
$\{1,1+1, 1+1+1 \}$, Brown et. al \cite{B} considered van der
Waerden's theorem on homothetic copies of
$\{1,1+s, 1+s+t \}$. 
 They defined $f(s,t)$ to be the smallest
positive integer $N$ such that every $2$-coloring of $[1, N]$ has
a monochromatic homothetic copy of $\{1,1+s, 1+s+t \}$. As they
have noted, $f(s,t)=f(t,s)$ and hence we assume that $s \geq t$.
They proved that $f(s,t) \leq 4(s+t)+1$ for all $s$, $t$ and that
the equality holds in the case where both $s/g\not\equiv0\pmod 4$
and $t/g\not\equiv0\pmod 4$ with $g=gcd(s,t)$
 and in
many other cases. Also they proved that $f(4mt,t)=4(4mt+t)-t+1$ or
$4(4mt+t)+1$. In this paper, we show that $f(4mt,t)=4(4mt+t)-t+1$
and that for all the other $(s,t)$, $f(s,t)=4(s+t)+1$.


The following theorem is stated as THEOREM 1 in Brown et. al
\cite{B}.

\begin{thm}
\label{thm:multiple} Let $s$, $t$ and $c$ be positive integers.
Then $f(cs, ct)=c(f(s,t)-1)+1$.
\end{thm}

\section{f(4mt, t)=4(4mt+t)-t+1}


\begin{lem}\label{a1}
 Let $N$, $s$, $t\in Z^+$ and $C:[1,N]\rightarrow\{0.1\}$ be a coloring
without 
 monochromatic triple $\{x,x+sy,x+(s+t)y\}$ for $y=1$, $2$ or $3$. Note that if $u \in \{0,1\}$ denotes one color of $C$, then $1-u$ denotes the other color.
 Define
$V(a)=(C(a),C(a+s),C(a+2s),C(a+3s))$ for $a \in [1, N-3s]$. Then
the followings are true for each $u\in\{0,1\}$.

\medskip\noindent
1) For $1\le a\le N-2s-2t$, if $C(a)=C(a+s)=u$, then
$C(a+2s)=1-u$.
\newline For $s+1\le a\le N-s-2t$, if
$C(a)=C(a+s)=u$, then $C(a-s)=1-u$.
\newline  For $2s+1\le a\le N-2t$, if $C(a)=C(a+t)=u$, then
$C(a+2t)=1-u$.
\newline For $2s+t+1\le a\le N-t$,
if $C(a)=C(a+t)=u$, then $C(a-t)=1-u$.

\medskip\noindent
 2) For $1\le a\le N-3s-4t$, if $V(a)=(u,1-u,1-u,u)$, then $V(a+t)=(1-u,1-u,u,u)$ or
$(1-u,u,u,1-u)$.

\medskip\noindent
 3) For $1\le a\le N-3s-6t$, if $V(a)=(u,u,1-u,u)$, then $V(a+t)=(1-u,1-u,u,u)$, $V(a+2t)=(1-u,u,u,1-u)$ and $V(a+3t)=(u,u,1-u,1-u)$.


\medskip\noindent
 4) For $1\le a\le N-3s-5t$, if $V(a)=(u,1-u,u,u)$,
then $V(a+2t)=(u,1-u,1-u,u)$ and $V(a+3t)=(1-u,u,u,1-u)$.

\medskip\noindent
 5) For $1\le a\le N-3s-4t$, if $V(a)=(u,1-u,1-u,u)$ and $V(a+t)=(1-u,1-u,u,u)$,
then $V(a+2t)=(1-u,u,u,1-u)$. Moreover, if $1\le a\le N-3s-6t$,
then $V(a+3t)=(u,u,1-u,1-u)$.
\newline For $1\le a\le N-3s-4t$,
if $V(a)=(u,1-u,1-u,u)$ and $V(a+t)=(1-u,u,u,1-u)$, then
$V(a+2t)=(u,1-u,1-u,u)$.
 \end{lem}

{\bf Proof.} Throughout the proof, we may assume that $u=0$.

\medskip\noindent
1) For $1\le a\le N-2s-2t$, if $C(a)=C(a+s)=0$, then $C(a+s+t)=1$.
Suppose $C(a+2s)=0$. Then
 $C(a+2s+t)=C(a+2s+2t)=1$.  This gives a monochromatic triple $\{a+s+t, a+2s+t, a+2s+2t\} \subset [1,N]$,
a contradiction. Hence $C(a+2s)=1$. The others are true by similar
methods.

\medskip\noindent
 2) For $1\le a\le N-3s-4t$, if $V(a)=(0,1,1,0)$, then $C(a+2s+t)=0$ and $C(a+3s+3t)=1$.
If $C(a+s+t)=0$,
 then $V(a+t)=(1,0,0,1)$ by 1).
If $C(a+s+t)=1$, then $C(a+3s+t)=0$ and hence $C(a+3s+2t)=1$. This
implies $C(a+3s+4t)=0$ and then $C(a+t)=1$. Thus
$V(a+t)=(1,1,0,0)$.

\medskip\noindent
 3) For $1\le a\le N-3s-6t$, if $V(a)=(0,0,1,0)$, then 
 $C(a+s+t)=C(a+3s+2t)=C(a+3s+3t)=1.$
Thus by 1), $C(a+3s+t)=C(a+3s+4t)=0.$ This implies $C(a+t)=1.$
Again by 1)
  $V(a+t)=(1,1,0,0)$. 
Now $C(a+s+2t)=0$. Also $C(a+3s+2t)=C(a+3s+3t)=1$ leads to
$C(a+2s+2t)=0$.
Again by 1) 
$V(a+2t)=(1,0,0,1)$.

\medskip\noindent
 4) For $1\le a\le N-3s-5t$, if $V(a)=(0,1,0,0)$, then
$C(a+2s+2t)=C(a+3s+t)=C(a+3s+3t)=1.$
 So $C(a+3s+2t)=C(a+s+t)=0.$
If $C(a+s+2t)=0$, then $C(a+t)=1$. And hence either $\{a+t,
a+3s+t, a+3s+4t\}$ or $\{a+s+2t, a+3s+2t, a+3s+4t\}$ is
monochromatic, a contradiction.
 Therefore
$C(a+s+2t)=1$ and hence $V(a+2t)=(0,1,1,0)$  
 by 1). Thus $C(a+2s+3t)=0$.
 $C(a+2t)=C(a+3s+2t)=0$ leads to $C(a+3s+5t)=1$
and hence it leads to $C(a+s+3t)=0.$
By 1), $V(a+3t)=(1,0,0,1)$.

\medskip\noindent
 5) Let $1\le a\le N-3s-4t$. Assume $V(a)=(0,1,1,0)$ and $V(a+t)=(1,1,0,0)$.
Then 
  $C(a+s+2t)=0$ 
and $C(a+3s+2t)=C(a+3s+3t)=1$. So $C(a+2s+2t)=0$ and hence
$V(a+2t)=(1,0,0,1)$.
If $1\le a\le N-3s-6t$, then 
$C(a+2s+3t)=1$ and $C(a+3s+4t)=C(a+3s+5t)=0$. So by 1)
 $C(a+3s+6t)=1$ and therefore $C(a+3t)= (0,0,1,1)$.
 Assume $V(a)=(0,1,1,0)$ and $V(a+t)=(1,0,0,1)$. Then $C(a+2s+2t)=C(a+3s+3t)=1$
and $C(a+3s+4t)=0$. So $C(a+3s+2t)=0$ and hence
 $C(a+s+2t)=1$.
Thus $V(a+2t)=(1,0,0,1)$.

\begin{lem}\label{b1}
Let $N$, $s$, $t$, $C$ and $V$ be the same as in Lemma \ref{a1}
and $a\in [1, N-3s]$. Then the followings are true for each
$u\in\{0,1\}$.

\medskip\noindent 1) If
$V(a)=(u,1-u,1-u,u)$ and $V(a+1)=(1-u,1-u,u,u)$, then
\begin{equation}\label{a2}
V(a+jt)=\left\{\begin{array}{ll}(u,1-u,1-u,u), & {j\equiv0\pmod4}
\\ (1-u,1-u,u,u), & {j\equiv1\pmod4} \\(1-u,u,u,1-u), & {j\equiv2\pmod4} \\
(u,u,1-u,1-u), &
{j\equiv3\pmod4}\\\end{array}\right.\end{equation} for all $j$
such that $0\le j\le2\lfloor{N-a-3s-2t\over2t} \rfloor.
$

\medskip\noindent
 2) If $V(a)=(u,1-u,1-u,u)$ and $V(a+1)=(1-u,u,u,1-u)$,
then \begin{equation}\label{a4}
V(a+jt)=\left\{\begin{array}{ll}(u,1-u,1-u,u), & {j\equiv0\pmod2}
\\ (1-u,u,u,1-u), & {j\equiv1\pmod2} \end{array}\right.\end{equation}
for all $j$ such that $0\le j\le\lfloor{N-a-3s-2t\over t}
\rfloor.$
\end{lem}

{\bf Proof.} Throughout the proof, we may assume that $u=0$.

\medskip\noindent
1) Let $k$ be the largest integer such that (\ref{a2}) holds for
all $0\le j\le k$. Assume $k\equiv0\pmod4$. Then
$V(k-2)=(1,0,0,1)$ and $V(k-1)=(0,0,1,1)$.
Suppose $k<2\lfloor{N-a-3s-2t\over2t}\rfloor$. Then 
 $a+(k-2)t\le N-3s-6t$
   and hence by 5) of Lemma \ref{a1} $V(k+1)=(1,1,0,0)$, a
contradiction. So $k\ge2\lfloor{N-a-3s-2t\over2t}\rfloor$.
Similarly, $k\ge2\lfloor{N-a-3s-2t\over2t}\rfloor$ if $k
\equiv1,2,3\pmod4$.
Thus (\ref{a2}) holds for all $0\le
j\le2\lfloor{N-a-3s-2t\over2t}\rfloor$.

\medskip\noindent
 2) Let $k$ be the largest integer such that (\ref{a4}) holds for all $0\le j\le
k$. Assume $k\equiv0\pmod2$. Then $V(k-1)=(1,0,0,1)$ and
$V(k)=(0,1,1,0)$. Suppose $k<\lfloor{N-a-3s-2t\over t}\rfloor$.
Then 
 $a+(k-1)t\le N-3s-4t$ and hence by 5) of Lemma \ref{a1}
$V(k+1)=(1,0,0,1)$, a contradiction. So
$k\ge\lfloor{N-a-3s-2t\over t}\rfloor$. Similarly,
$k\ge\lfloor{N-a-3s-2t\over t}\rfloor$ if $k\equiv1\pmod2$.
 Thus
(\ref{a4}) holds for all $0\le j\le\lfloor{N-a-3s-2t\over
t}\rfloor$.


\begin{thm}
\label{thm: f(4m,1)} For all $m\in Z^+$, $f(4m,1)=16m+4$.
\end{thm}

{\bf Proof.} THEOREM $4$ of Brown et. al \cite{B} states that
$f(4mt,t)=4(4mt+t)-t+1$ or $4(4mt+t)+1$ for all positive integers
$m$ and $ t$.
So it suffices to show that $f(4m,1)<16m+5$. Suppose not. Then
there is a coloring of [1, 16m+4] which has no monochromatic
homothetic copy of $\{1, 4m+1, 4m+2\}$. Let $C$ be such a coloring
and $V(a)=(C(a), C(a+4m),C(a+8m),C(a+12m))$ for
 each $a \in [1, 4m+4]$.
 We may assume that $C(1)=0$.
There are five cases to consider. In each case we show that we
have contradictions.

\medskip\noindent
Case 1: $V(1)=(0,0,1,0)$. Then by 3) of Lemma
\ref{a1} $V(2)=(1,1,0,0)$, $V(3)=(1,0,0,1)$ and $V(4)=(0,0,1,1)$.
So by 1) of Lemma \ref{b1}
\begin{equation}V(3+j)=\left\{\begin{array}{ll}(1,0,0,1), &
{j\equiv0\pmod4}
\\ (0,0,1,1), & {j\equiv1\pmod4} \\(0,1,1,0), & {j\equiv2\pmod4} \\
(1,1,0,0), & {j\equiv3\pmod4}\\
\end{array}\right.\end{equation}
for all $j\le4m-2$. In particular, $V(4m+1)=(0,1,1,0)$, which
contradicts to $C(12m+1)=0$.

\medskip\noindent
Case 2: $V(1)=(0,1,0,0)$. Then by 4) of Lemma \ref{a1}
$V(3)=(0,1,1,0)$ and $V(4)=(1,0,0,1)$. So by 2) of Lemma \ref{b1}
 \begin{equation}
V(3+j)=\left\{\begin{array}{ll}(0,1,1,0), & {j\equiv0\pmod2}
 \\(1,0,0,1), & {j\equiv1\pmod2} \\\end{array}\right.\end{equation}
for all $j\le4m-1$. In particular,  $V(4m+1)=(0,1,1,0)$, which
contradicts to $C(12m+1)=0$.

\medskip\noindent
Case 3: $V(1)=(0,1,1,0)$. Then by 2) of Lemma \ref{a1} $V(2)$ is
$(1,1,0,0)$ or $(1,0,0,1)$. If $V(2)=(1,1,0,0)$, then by 1) of
LEMMA \ref{b1}
\begin{equation}V(1+j)=\left\{\begin{array}{ll}(0,1,1,0), &
{j\equiv0\pmod4}\\ (1,1,0,0), & {j\equiv1\pmod4} \\(1,0,0,1), & {j\equiv2\pmod4} \\
(0,0,1,1), & {j\equiv3\pmod4}\\\end{array}\right.\end{equation}
for all $j\le4m$. In particular, $V(4m+1)=(0,1,1,0)$, which
contradicts to $C(12m+1)=0$.
\newline If $V(2)=(1,0,0,1)$, then by
2) of Lemma \ref{b1}  \begin{equation}
V(1+j)=\left\{\begin{array}{ll}(0,1,1,0), & {j\equiv0\pmod2}
 \\(1,0,0,1), & {j\equiv1\pmod2} \\\end{array}\right.\end{equation}
for all $j\le4m+1$. In particular, $V(4m+1)=(0,1,1,0)$, which
contradicts to $C(12m+1)=0$.

\medskip\noindent
Case 4: $V(1)=(0,0,1,1)$. Then $V(4m+1)\not=(0,0,1,1)$ since
$C(8m+1)=1$. So there is a minimal $k$ such that
$V(k)\not=(u(k),u(k),1-u(k),1-u(k))$ where $u(k)=
 \begin{cases}0, &k \text{ is odd} \\
1, &k\text{ is even} \end{cases}$.

 Consider the case when $k$ is even.
As $V(k-1)=(0,0,1,1)$, $V(k)$ is
$(0,1,0,0)$ or $(0,1,1,0)$. If $k=4m$, then it
contradicts to $C(12m+2)=0$ and $C(12m+1)=1$.
 Assume $k\le4m-2$.

Firstly if $V(k)=(0,1,0,0)$, then by 4) of Lemma \ref{a1} $V(k+2)=
(0,1,1,0)$ and $V(k+3)=(1,0,0,1)$. So by  2) of Lemma \ref{b1}
\begin{equation}
V(k+2+j)=\left\{\begin{array}{ll}(0,1,1,0), & {j\equiv0\pmod2}
 \\(1,0,0,1), & {j\equiv1\pmod2} \\\end{array}\right.\end{equation}
for all $j\le4m-k$. In particular, $V(4m+1)=(1,0,0,1)$, which
contradicts to $C(12m+1)=1$.

Secondly if $V(k)=(0,1,1,0)$, then by 2) of Lemma \ref{a1}
$V(k+1)$ is $(1,1,0,0)$ or $(1,0,0,1)$. If $V(k+1)=(1,1,0,0)$,
then by 1) of Lemma \ref{b1}
\begin{equation}V(k+j)=\left\{\begin{array}{ll}(0,1,1,0), &
{j\equiv0\pmod4}\\(1,1,0,0), & {j\equiv1\pmod4} \\(1,0,0,1), & {j\equiv2\pmod4} \\
 (0,0,1,1), & {j\equiv3\pmod4}\\\end{array}\right.\end{equation}
for all $j\le4m-k+2$. So  $V(4m+1)$
is $(1,1,0,0)$ or $(0,0,1,1)$, which contradicts to $C(12m+1)=C(8m+1)=1$.
 If $V(k+1)=(1,0,0,1)$, then by 2) of Lemma \ref{b1}
 \begin{equation} V(k+j)=\left\{\begin{array}{ll}(0,1,1,0), &
{j\equiv0\pmod2}
 \\(1,0,0,1), & {j\equiv1\pmod2} \\\end{array}\right.\end{equation}
for all $j\le4m-k+2$. In particular, $V(4m+1)$ is $(1,0,0,1)$,
which contradicts to $C(8m+1)=1$.

 When $k$ is odd,
then we get contradictions similarly.
%
%
%

\medskip\noindent
 Case 5: $V(1)=(0,1,0,1)$. Then $C(8m+3)=1$ and $C(12m+3)=0$. Thus
$V(2)$ is $(0,0,1,0)$, $(0,0,1,1)$, $(0,1,0,1)$, $(1,0,1,1)$,
$(1,1,0,1)$, $(1,0,1,0)$ or $(1,0,0,1)$.

If $V(2)=(0,0,1,0)$, then by 3) of Lemma \ref{a1}
$V(3)=(1,1,0,0)$, which contradicts to $C(8m+3)=1$.

If $V(2)=(1,1,0,1)$, then similarly $V(3)=(0,0,1,1)$, which contradicts to $C(12m+3)=0$.

If $V(2)=(0,0,1,1)$,
then 
$C(4m+3)=1$ and hence by 1) of Lemma \ref{a1}
  $V(3)=(0,1,1,0)$. By the method used in Case 4, we
can show that $V(4m+1)$ is $(0,1,1,0)$ or $(1,0,0,1)$, which
contradicts to $C(4m+1)=C(12m+1)=1$.

If $V(2)=(1,0,1,1)$, then by 4) of Lemma \ref{a1} $V(4)=(1,0,0,1)$
and $V(5)=(0,1,1,0)$. So by 2) of Lemma \ref{b1}
\begin{equation} V(4+j)=\left\{\begin{array}{ll}(1,0,0,1), &
{j\equiv0\pmod2}
\\(0,1,1,0), & {j\equiv1\pmod2}\\\end{array}\right.\end{equation}
for all $j\le4m-2$. So $V(4m+1)=(0,1,1,0)$, which contradicts to
$C(4m+1)=1$.

If $V(2)=(1,0,0,1)$, then by 1) and 2) of Lemma \ref{a1} 
$V(3)=(0,1,1,0)$. So by 2) of Lemma \ref{b1}
 \begin{equation}V(2+j)=\left\{\begin{array}{ll}(1,0,0,1), & {j\equiv0\pmod2}
 \\(0,1,1,0), & {j\equiv1\pmod2} \\\end{array}\right.\end{equation}
for all $j\le4m$. In particular, $V(4m+1)=(0,1,1,0)$, which
contradicts to $C(4m+1)=1$.

If
 $V(2)=(0,1,0,1)$, then
$C(8m+3)=C(8m+4)=1$ and $C(12m+3)=C(12m+4)=0$ and hence
$C(4m+3)=0$. $C(3)=1$ as otherwise $V(3)=(0,0,1,0)$ and hence by
3) of Lemma \ref{a1} $V(4)=(1,1,0,0)$, a contradiction. Thus
$V(3)=(1,0,1,0)$. As $C(4m+1)=1$, $V(4m+1)\not=(0,1,0,1).$ So
there is a minimal $k$ such that
$V(k)\not=(u(k),1-u(k),u(k),1-u(k))$ where $u(k)=\begin{cases}0,
&k\equiv1,2\pmod4 \\1, &k\equiv0,3\pmod4\end{cases}$.
Firstly assume $k\equiv1\pmod4 $. Then $k\le 4m+1$.
When $k\le4m-3$, $V(k-2)=V(k-1)=(1,0,1,0)$ and hence as above 
$V(k)=(1,1,0,1)$, a contradiction. When $k=4m+1$,
 $C(8m+k=12m+1)=0$ as $V(k-2)=(1,0,1,0)$, a
contradiction.
Secondly assume $k\equiv2\pmod4$. Then $k\le4m-2$. From
$V(k-2)=(1,0,1,0)$ and $V(k-1)=(0,1,0,1)$, we get
$C(8m+k)=C(12m+k+1)=0$ and $C(8m+k+1)=C(12m+k)=1$.
 Suppose $C(4m+k)=1$. Then $C(12m+k+2)=0$ and hence by 1) of Lemma \ref{a1} $C(12m+k+3)=1$.
 Thus $C(k)=0$ i.e. $V(k)=(0,1,0,1)$,
a contradiction. Therefore $C(4m+k)=0$ and hence $V(k)=(1,0,0,1)$
by 1) of lemma \ref{a1}.
By 2) of Lemma \ref{a1} 
 $V(k+1)=(0,1,1,0).$ So by
2) of Lemma \ref{b1}
\begin{equation}V(k+j)=\left\{\begin{array}{ll}(1,0,0,1), & {j\equiv0\pmod2}
 \\(0,1,1,0), & {j\equiv1\pmod2} \\\end{array}\right.\end{equation}
for all $j\le4m-k+2$. In particular, $V(4m+1)=(0,1,1,0)$, which
contradicts to $C(4m+1)=1$.
Thirdly assume $k\equiv3\pmod4$. Then
$k\le4m-1.$  
When $k\le4m-5$, as in the case when $k\equiv1\pmod4$, we get a
contradiction.
When $k=4m-1$, $V(k+1)=(0,1,1,0)$ and hence $C(8m+k+2=12m+1)=0$, a
contradiction.
Lastly
 assume $k\equiv0\pmod4$. Then $k\le4m-4$.
As in the case when $k\equiv2\pmod4$, we get contradictions.

If $V(2)=(1,0,1,0)$, then $V(4m+1)\not=(0,1,0,1)$ since $C(4m+1)=1$.
So there is a minimal $k$ such that $V(k)\not=(u(k),1-u(k),u(k),1-u(k))$ where
$u(k)=\begin{cases}0, &k\equiv0,1\pmod4 \\
1, &k\equiv2,3\pmod4\end{cases}$. For each $k$, the case reduces
to the previous case where $V(2)=(0,1,0,1)$ for $k+1$ and we get
contradictions.

\begin{rmk} 
 Investigating the proof of Theorem \ref{thm:
f(4m,1)}, we see that we may have a coloring $C$ of $[1, 16m+3]$
such that $C(1)=0$ which has no monochromatic homothetic copy of
$\{1, 4m+1, 4m+2\}$ only in Case 4, where either $k$ is even,
$k\le4m-2$, $V(k)=(0,1,1,0)$ and $V(k+1)=(1,1,0,0)$ or $k$ is odd,
$k\le4m-1$, $V(k)=(1,0,0,1)$ and $V(k+1)=(0,0,1,1)$. In fact only
in the latter, if
$k\equiv3\pmod4$, 
we have such a coloring which satisfies
 \begin{equation}V(j)=\left\{\begin{array}{ll}(0,0,1,1), & {j\equiv1\pmod2, j<k}
 \\(1,1,0,0), & {j\equiv0\pmod2,j<k} \\(1,0,0,1), & {j\equiv3\pmod4,k\le j\le 4m}
 \\(0,0,1,1), & {j\equiv0\pmod4,k\le j\le 4m} \\(0,1,1,0), & {j\equiv1\pmod4,k\le j\le 4m}
 \\(1,1,0,0), & {j\equiv2\pmod4,k\le j\le 4m}\\\end{array}\right.\end{equation},
 $C(16m+1)=C(16m+2)=0$ and $C(16m+3)=1$.
We can dispose it as
\newline$\begin{array}{ccccc}
  0 & 0 & 1 & 1 &0 \\1 & 1 & 0 & 0 &0 \\ 0 & 0 & 1 & 1 &1 \\ 1 & 1 & 0 & 0 & \\
  & \vdots & \vdots & & \\ 0 & 0 & 1 & 1 & \\1 & 1 & 0 & 0& \\1& 0 & 0 & 1& \\
0 & 0 & 1 & 1 & \\ 0 & 1 & 1 & 0 & \\1 & 1 & 0 & 0&\\  & \vdots &\vdots & & \\
 1& 0 & 0 &1 & \\ 0 & 0 & 1 & 1& \end{array}$
\newline
 with
   $C(4m(j-1)+i)$ in the $i$-th row and the $j$-th
column.
 For each $C$, $1-C$ has no monochromatic homothetic copy of $\{1,4m+1, 4m+2\}$.
 Thus
we have altogether $2m$ colorings of $[1, 16m+3]$
which has no monochromatic homothetic copy of $\{1,4m+1, 4m+2\}$.


\end{rmk}

 Thus together with Theorem \ref{thm:multiple}, the following theorem is true.

\begin{thm}
\label{thm:f(4mt,t)} For all $m$, $t\in Z^+$,
$f(4mt,t)=4(4mt+t)-t+1$.
\end{thm}

\section{Remaining cases}

The following theorem which is stated as
 THEOREM $3$ in Brown et. al\cite{B} implies that
 $f(s,t)=4(s+t)+1$ if $t$ divides $s$ and $s/t \not\equiv0\pmod 4$.

\begin{thm}
\label{thm:THEOREM 3} Let $s$, $t$ be positive integers with
$g=gcd(s,t)$. If $s/g\not\equiv0\pmod 4$ and $t/g\not\equiv0\pmod
4$, then $f(s,t)=4(s+t)+1$.
\end{thm}


\begin{thm}
\label{thm:s/t} Let $s$, $t$ be positive integers such that $s
> t > 1$ and $t$ does not divide $s$. If $\lfloor s/t \rfloor$ is even or
$\lfloor 2s/t \rfloor$ is even, then $f(s,t)=4(s+t)+1$. If
$\lfloor s/t \rfloor$ and $\lfloor 2s/t \rfloor$ are both odd,
then $f(s,t)=4(s+t)+1$ provided $s$, $t$ satisfy the additional
condition
 $s/t \notin (1.5,2)$.
\end{thm}

This theorem is stated as THEOREM $5$ in Brown et. al \cite{B}.
Together with the following theorem it implies that
 $f(s,t)=4(s+t)+1$ if $t$ does not divide $s$.

\begin{thm}
\label{thm:s/t in} Let $s$, $t$ be positive integers such that $s
> t > 1$ and $t$ does not divide $s$. If $\lfloor s/t \rfloor$ and
$\lfloor 2s/t \rfloor$ are both odd and $s/t \in (1.5,2)$, then
$f(s,t)=4(s+t)+1$.
\end{thm}

The following lemma is directly from the Chinese Remainder
Theorem.

\begin{lem} \label{lem:i,j} Let $s$, $t$ be positive integers such
that $gcd(s, t) = 1$. Then any integer in $[1, 4s+4t] $ is equal
to $1+is+jt$ for some unique $0 \leq j < s$ and $-jt/s \leq i <
4+(4-j)t/s$.
\end{lem}

{\bf Proof of Theorem \ref{thm:s/t in}.} By Theorem
\ref{thm:multiple}, it is enough to consider $s$, $t$ such that
$gcd(s,t)=1$. By Theorem \ref{thm:THEOREM 3}, we need to consider
only the cases where $s\equiv0\pmod 4$ or $t \equiv0\pmod 4$. By
THEOREM $2$ of Brown et. al \cite{B}, $f(s,t) \leq 4(s+t)+1$ for
all positive integers $s$ and $t$.
In each of the following exhaustive cases, we show that the
equality holds by
 choosing a $2-$coloring $C$ of $[1, 4s+4t]$ which contains no
homothetic copy of $\{1,1+s, 1+s+t \}$.
By Lemma \ref{lem:i,j}, for each $1+is+jt \in [1, 4s+4t]$, $ 0
\leq j < s$ and $-jt/s \leq i < 4+(4-j)t/s$. Figure $1$-Figure $5$
show $C(1+is+jt)$
 for those $i$ and $j$.

\medskip\noindent
Case 1: $s \equiv0\pmod 4$ and $t \equiv3\pmod 4$. Consider the
following two subcases.

\medskip\noindent Case 1a: $s/t \geq 5/3$. Note that
$1+(-t+7)s+(s-1)t \geq 1+4s+4t$ and hence $1+(-t+7)s+(s-1)t \notin
[1,4s+4t]$. We marked $@$ on its position in Figure 1.
 Let $A=\{1+is+jt \mid i\geq -t+7$ or $j\leq s-4\}$ and
$B=\{1+is+jt \mid i\leq -t+6, j \geq s-3\}$.
For $1+is+jt \in A$ let

\smallskip\noindent
 $C(1+is+jt)=$ $\left\{\begin{array}{ll}
    0, & \hbox{$j$ even, $i\equiv0,1\pmod 4;$ $j$ odd, $i \equiv2,3\pmod 4$} \\
    1, & \hbox{otherwise} \\
\end{array}
\right.$

\smallskip\noindent
 and for $1+is+jt \in B$ let

\smallskip\noindent $C(1+is+jt)=$
$\left\{%
\begin{array}{ll}
    0, & \hbox{$j$ even $i\equiv1,2\pmod 4$; $j$ odd $i$ $\equiv0,3\pmod 4$} \\
    1, & \hbox{otherwise}
    \\
\end{array}%
\right.$.

 To show that $C$ has no monochromatic homothetic copy of $\{1,
1+s, 1+s+t\}$, it is enough to show that $C$ has no monochromatic
triple $\{1+is+jt, 1+(i+y)s+jt, 1+(i+y)s+(j+y)t\} $ where $y=1$,
$2$ or $3$. Suppose $C$ has such a triple. Then the triple is not
in $A$ as $C(1+(i+y)s+jt) \neq C(1+(i+y)s+(j+y)t)$ for $y=1$, $3$
and $C(1+(i+y)s+jt) \neq C(1+is+jt)$ for $y=2$. Similarly the
triple is not in $B$ either.
Note that $1+(i+y)s+jt \in A$ from the definitions of $A$ and $B$.
If $y=1$, then $1+(i+y)s+(j+y)t \in B$ and hence $i \leq -t+5$ and
$j = s-4$. As $i\geq -jt/s$, $i \geq
-t+3$. 
For each of such an $(i,j)$, $\{1+is+jt, 1+(i+y)s+jt,
1+(i+y)s+(j+y)t\}$ is not monochromatic by the definitions of $A$
and $B$, a contradiction. By a similar argument, the triple is not
monochromatic if $y=3$.
If $y=2$, then $1+is+jt \in B$. Thus $i \geq -t+5$ and $j \geq
s-3$ and therefore $1+(i+2)s+(j+2)t \geq
1+4s+4t$, a contradiction. 
Therefore $C$ avoids monochromatic homothetic copy of $\{1, 1+s,
1+s+t\}$.

\begin{figure}
  \includegraphics[width=6.0in]{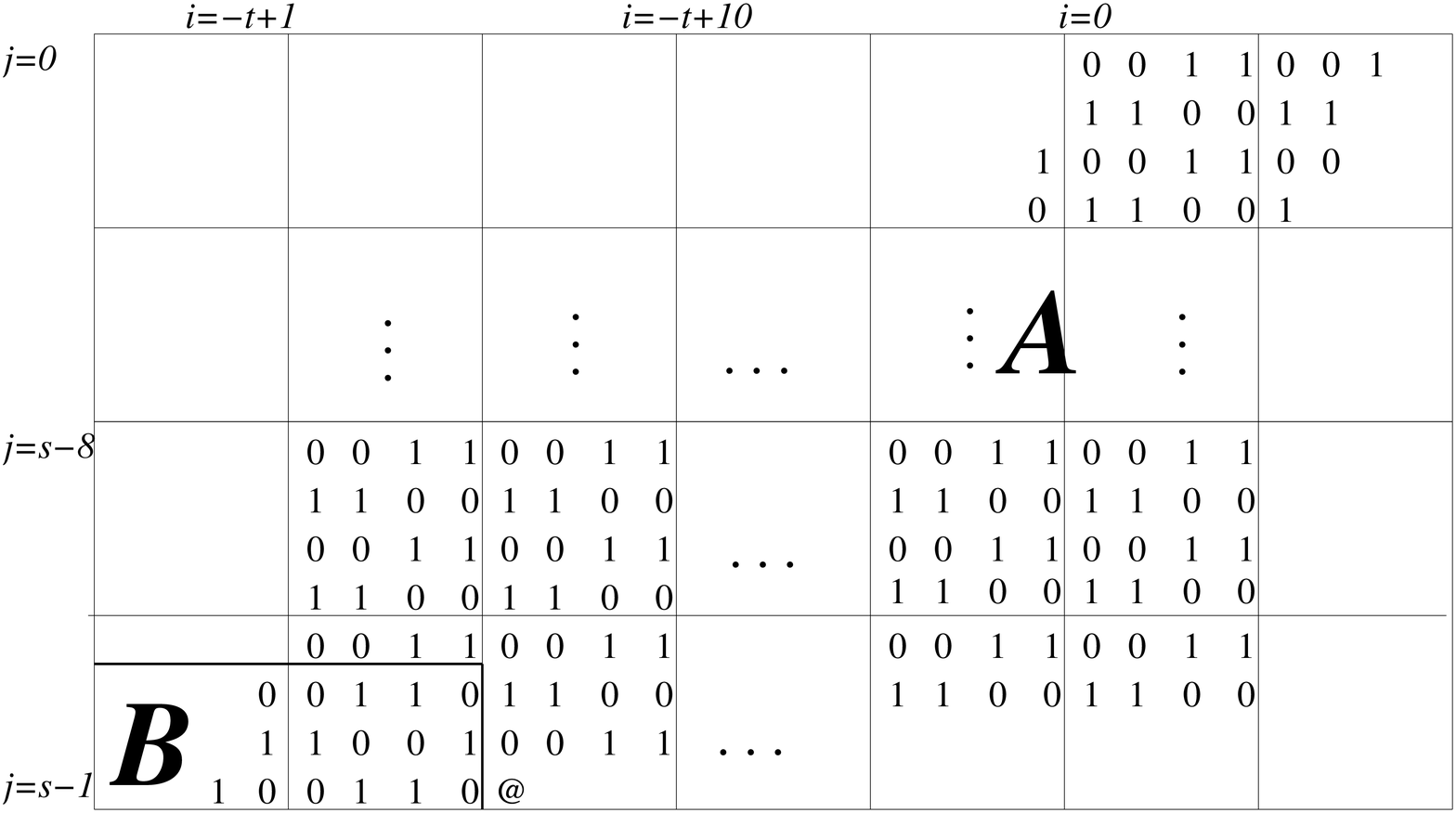}\\
  \caption{ Case 1a}
\end{figure}

The following $2-$coloring $C$ of $[1, 4s+4t]$
 in each case has no monochromatic homothetic copy of $\{1, 1+s, 1+s+t\}$ by similar arguments.

\medskip\noindent
 Case 1b: $s/t < 5/3$. Note that $1+(-t+3)s+(s-5)t < 1$ and hence $1+(-t+3)s+(s-5)t \notin [1, 4s+4t]$.
We marked $@$ on its position in Figure 2.
Let $C(1+4s+2t)=1$ and $C(1+5s+2t)=C(1+4s+3t)=0$. For the other
$(i, j)$ if $i \geq -t+7$ or $2 \leq j \leq s-5$, then let

\noindent $C(1+is+jt)=$
$\left\{%
\begin{array}{ll}
    0, & \hbox{$j$ even $i \equiv0,3\pmod 4$; $j$ odd $i \equiv1,2\pmod 4$} \\
    1, & \hbox{otherwise} \\
\end{array}%
\right.$.

\noindent
For the other $(i, j)$, let

\noindent $C(1+is+jt)=$
$\left\{%
\begin{array}{ll}
    0, & \hbox{$(i,j)=(0,0)$, $(2,0)$, $(5,0)$, $(0,1)$, $(1,1)$, $(4,1)$, $(-t+4, s-4)$},\\
& \hbox{ $(-t+3,s-3)$, $(-t+5,s-3)$, $(-t+6,s-3)$, $(-t+2,s-2)$}, \\
    & \hbox{ $(-t+4,s-2)$, $(-t+5, s-2)$,$(-t+1,s-1)$, $(-t+3,s-1)$},\\
    & \hbox{ $(-t+4,s-1)$, $(-t+6, s-1)$},\\
    1, & \hbox{otherwise}\\
\end{array}%
\right.$.

\begin{figure}
  \includegraphics[width=6.0in]{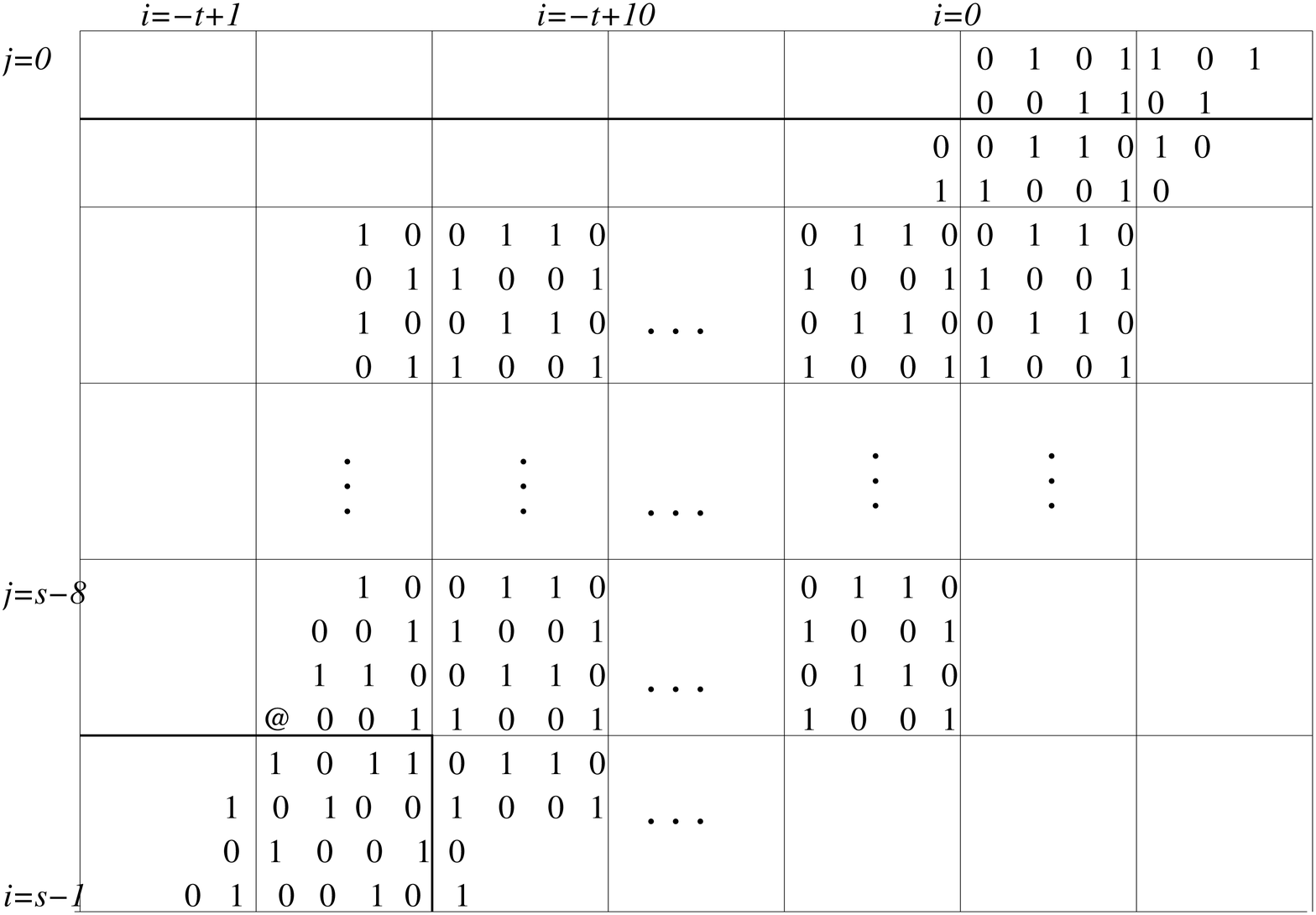}\\
  \caption{ Case 1b}
\end{figure}

\medskip\noindent
Case 2: $s \equiv0\pmod 4$ and $t \equiv1\pmod 4$. Consider the
following two subcases.

\medskip\noindent
Case 2a: $s/t \geq 5/3$. As in Case 1a, $1+(-t+7)s+(s-1)t \notin
[1, 4s+4t]$. We marked $@$ on its position in Figure 3.
Let $C(1+3s+4t)=1$, $C(1+3s+5t)=0$ and $C(1+6s)=0$. For the other
$(i, j)$, if $4 \leq j \leq s-4$, then let

\noindent $C(1+is+jt)=$
$\left\{%
\begin{array}{ll}
    0, & \hbox{$j$ even $i \equiv2,3\pmod 4$; $j$ odd $i\equiv0,1\pmod 4$} \\
    1, & \hbox{otherwise} \\
\end{array}%
\right.$.

\noindent
For the other $(i,j)$, if $ j \leq 3$, then let

\noindent $C(1+is+jt)=$
$\left\{%
\begin{array}{ll}  0, & j
    \hbox{ even $i \equiv1,2\pmod 4$; $j$ odd $i \equiv0,3\pmod 4$} \\
    1, & \hbox{otherwise} \\
\end{array}%
\right.$.

\noindent
For the other $(i,j)$, let

\noindent $C(1+is+jt)=$
$\left\{%
\begin{array}{ll}
    0, & \hbox{$(i,j)=(-t+2, s-3)$, $(-t+3,s-3)$, $(-t+6,s-3)$}, \\
    & \hbox{$(-t+4,s-2)$, $(-t+5,s-2)$, $(-t+7, s-2)$, $(-t+3,s-1)$},    \\
    & \hbox{ $(-t+4,s-1)$, $(-t+6, s-1)$,} \\
    1, & \hbox{otherwise} \\
\end{array}%
\right.$.

\begin{figure}
  \includegraphics[width=6.0in]{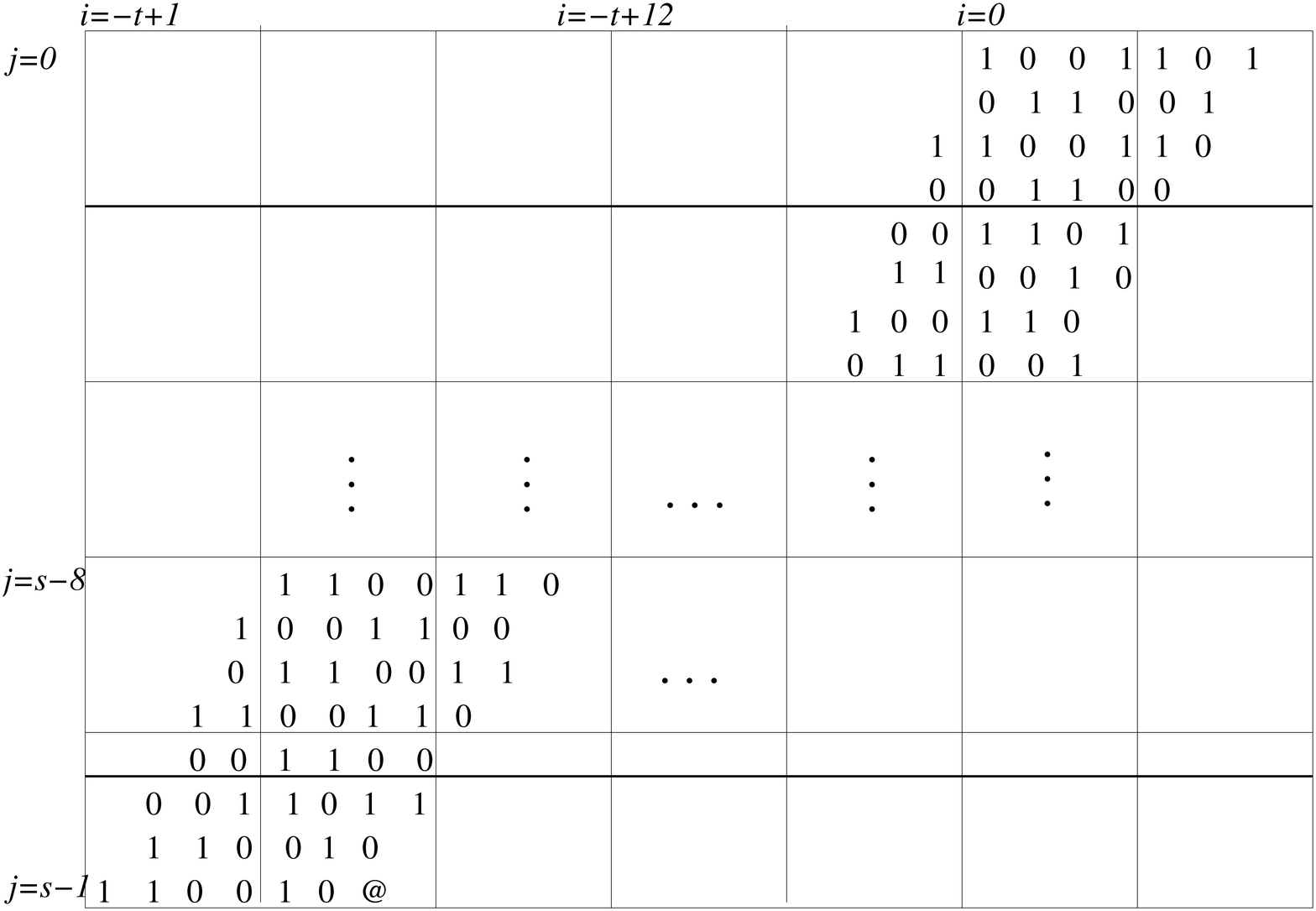}\\
  \caption{ Case 2a}
\end{figure}

\medskip\noindent
Case 2b: $s/t < 5/3$. As in Case 1b, $1+(-t+3)s+(s-5)t \notin [1,
4s+4t] $. We marked $@$ on its position in Figure 4.
Let $C(1+(-t+4)s+(s-6)t)=1$, $C(1+(-t+4)s+(s-5)t)=0$,
$C(1+(-t+8)s+(s-4)t)=0$ and $C(1+(-t+8)s+(s-3)t)=1$. For the other
$(i,j)$, if $i \leq 3$ and $2 \leq j \leq s-5$, then let

\noindent $C(1+is+jt)=$
$\left\{%
\begin{array}{ll}
    0, & \hbox{$j$ even $i\equiv0,3\pmod 4$; $j$ odd $i \equiv1,2\pmod 4$} \\
    1, & \hbox{otherwise} \\
\end{array}%
\right.$.

\noindent
For the other $(i,j)$, if $ j \geq s-4$, then let

\noindent $C(1+is+jt)=$
$\left\{%
\begin{array}{ll}
    0, & \hbox{$j$ even $i \equiv0,1\pmod 4$; $j$ odd $i \equiv2,3\pmod 4$} \\
    1, & \hbox{otherwise} \\
\end{array}%
\right.$.

\noindent
For the other $(i,j)$, let

\noindent $C(1+is+jt)=$
$\left\{%
\begin{array}{ll}
    0, & \hbox{$(i,j)=(0,0)$, $(2,0)$, $(5,0)$, $(6,0)$, $(1,1)$, $(3,1)$,
$(4,1)$, $(5,2)$, $(4,3)$} ,\\
    1, & \hbox{otherwise} \\
\end{array}%
\right.$.

\begin{figure}
  \includegraphics[width=6.0in]{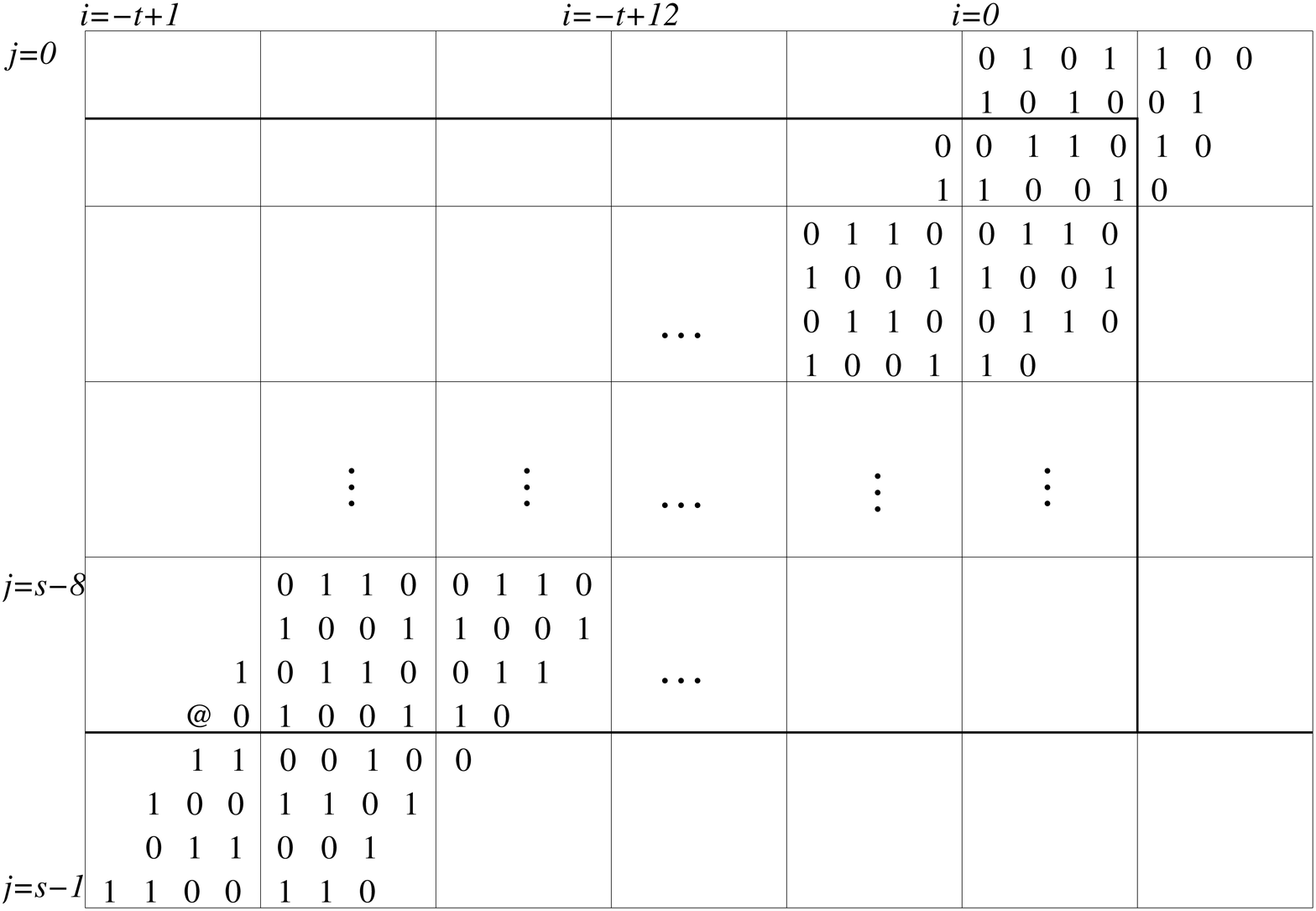}\\
  \caption{ Case 2b}
\end{figure}

\medskip\noindent
Case 3: $s \equiv1, 3\pmod 4$ and $t \equiv0\pmod 4$. Notice that
either $1+(-t+7)s+(s-1)t\notin [1, 4s+4t]$ or $1+(-t+3)s+(s-5)t
\notin [1, 4s+4t]$ depends on $s/t \geq 5/3$ or not. We marked
$@'$ and $@$ on their positions
 in Figure 5.
Let $C(1+3s+4t)=1$, $C(1+3s+5t)=0$, $C(1+6s)=1$, $C(1+s-t)=1$,
$C(1+2s-t)=1$, $C(1+3s-t)=0$, $C(1+4s-t)=0$, $C(1+5s-t)=1$,
$C(1+6s-t)=0$, $C(1+7s-t)=0$. For the other $(i,j)$, if $4 \leq j
\leq s-2$, then let

\noindent $C(1+is+jt)=$
$\left\{%
\begin{array}{ll}
    0, & \hbox{$j$ even $i \equiv2,3\pmod 4$; $j$ odd $i \equiv0,1\pmod 4$} \\
    1, & \hbox{otherwise} \\
\end{array}%
\right.$.

\noindent
For the other $(i,j)$, if $ j \leq 3$, then
 let

\noindent
 $C(1+is+jt)=$
$\left\{%
\begin{array}{ll}
    0, & \hbox{$j$ even $i \equiv1,2\pmod 4$; $j$ odd $i\equiv0,3\pmod 4$} \\
    1, & \hbox{otherwise} \\
\end{array}%
\right.$.


\begin{figure}
  \includegraphics[width=6.0in]{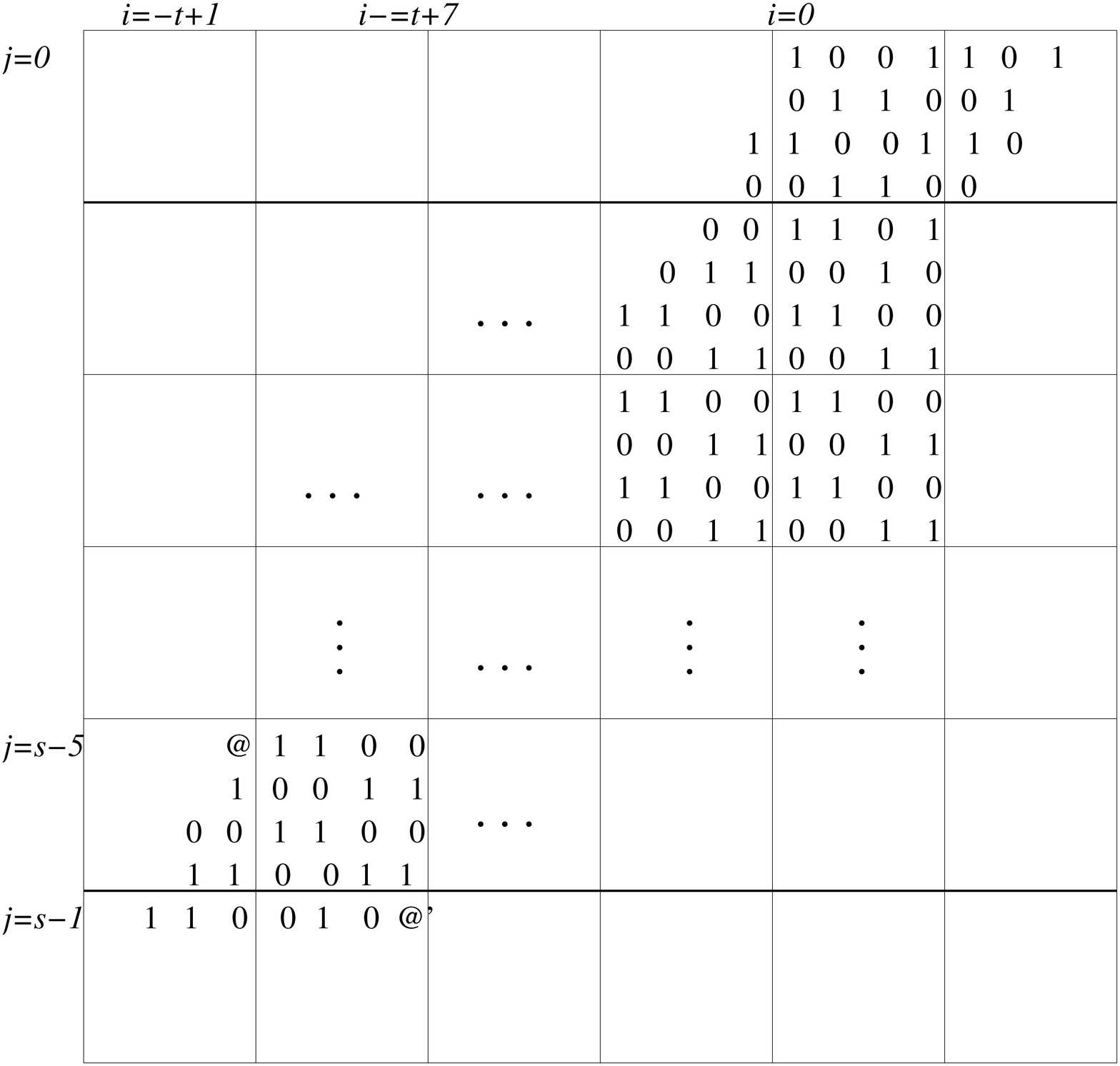}\\
  \caption{ Case 3}
\end{figure}

By Theorem \ref{thm:THEOREM 3}, Theorem \ref{thm:s/t} and Theorem
\ref{thm:s/t in}, for all pair of positive integers $(s,t)\neq
(4mt,t)$ for a positive integer $m$, $f(s,t)=4(s+t)+1$.

{\it Acknowledgement: We thank Seoul National University for
providing us with a place to research during the winter break of
2004.}


\begin{thebibliography}{4}



\bibitem{B};
T. C. Brown, B. M. Landman and M. Mishna,
\newblock  {\it Monochromatic homothetic copies of $\{1,1+s,1+s+t\}$},
\newblock  Canad. Math. Bull., 40. No.2 (1997); 149-157.

\bibitem{W};
B. L. Van der Waerden,
\newblock {\it Beweis einer Baudetschen Vermutung},
\newblock Nieuw Arch. Wisk, 15 (1927); 212-216.

\bibitem{GRS};
R. L. Graham, B. L. Rothschild and J. H. Spencer,
\newblock {\it Ramsey Theory},
\newblock  Wiley-Interscience, New York, 1990.



\end{thebibliography}
\end{document}